\documentclass{article}

\usepackage[preprint]{neurips_2023}

\usepackage{graphicx}

\usepackage[utf8]{inputenc} 
\usepackage[T1]{fontenc}    
\usepackage{hyperref}       
\usepackage{url}            
\usepackage{booktabs}       
\usepackage{amsfonts}       
\usepackage{nicefrac}       
\usepackage{microtype}      
\usepackage{xcolor}         
\usepackage{amsmath}        
\usepackage{mathtools}        
\usepackage{amssymb}        
\usepackage{algorithm}
\usepackage[noend]{algpseudocode}

\title{Constrained Bayesian optimization with merit functions}

\author{%
   Jingyi Wang \\
  Lawrence Livermore National Laboratory\\
  Livermore, CA 94550 \\
  \texttt{wang125@llnl.gov} \\
   \AND
   Cosmin G. Petra \\
  Lawrence Livermore National Laboratory\\
  Livermore, CA 94550 \\
   \texttt{petra1@llnl.gov} \\
   \And
   J. Luc Peterson \\
  Lawrence Livermore National Laboratory\\
  Livermore, CA 94550 \\
   \texttt{peterson76@llnl.gov} \\
}

\begin{document}

\newcommand{\Rbb}{\ensuremath{\mathbb{R} }}
\newcommand{\Ebb}{\ensuremath{\mathbb{E} }}
\newcommand{\norm}[1]{\left\lVert {#1} \right\rVert}

\maketitle

\begin{abstract}
  Bayesian optimization is a powerful optimization tool for problems where native first-order derivatives are unavailable. 
  Recently, constrained Bayesian optimization (CBO) has been applied to many engineering applications where constraints are essential.
  However, several obstacles remain with current CBO algorithms that could prevent a wider adoption. 
  We propose CBO algorithms using merit functions, such as the penalty merit function, in acquisition functions,
  inspired by nonlinear optimization methods, \textit{e.g.}, sequential quadratic programming.
  Merit functions measure the potential progress of both the objective and constraint functions, thus increasing algorithmic efficiency and allowing infeasible initial samples. 
  The acquisition functions with  merit functions are relaxed to have closed forms, making its implementation readily available wherever Bayesian optimization is. 
  We further propose a unified CBO algorithm that can be seen as extension to the popular expected constrained improvement (ECI) approach. 
We demonstrate the effectiveness and efficiency of the proposed algorithms through numerical experiments on synthetic problems and a practical data-driven engineering design problem in the field of plasma physics.
\end{abstract}

\section{Introduction}
Bayesian optimization is a black-box optimization method that does not require derivative information of the objective~\cite{bosurvey2023}.
Instead, a surrogate model, often a Gaussian process (GP), is constructed to replace the objective function based on iterative sampling of the data.  
It has been widely adopted in many applications such as inverse problems~\cite{wang2004}, structural design~\cite{mathern2021}, additive manufacturing process simulation~\cite{wang2023optimization}, robotics~\cite{calandra2016}, inertial confinement fusion (ICF) design~\cite{wang2023multifidelity}, etc. 
The method has been successfully extended to multifidelity Bayesian optimization with multifidelity surrogate models~\cite{zuluaga2013} and Gaussian processes with independent constraints~\cite{bernardo2011}.

Constrained Bayesian optimization (CBO) has attracted huge interest recently since constraints are critical to many engineering applications. The mathematical formulation of an inequality constrained optimization problem is  
\begin{equation} \label{eqn:opt-prob}
 \centering
  \begin{aligned}
	  &\underset{\substack{x}\in C}{\text{minimize}} 
	  & & f(x), \\
          &\text{subject to}
	  & & c(x) \geq  0,
  \end{aligned}
\end{equation}
where $C\subset \Rbb^n$ represents the bound constraints on $x$, $f:\Rbb^n\to \Rbb$ is the objective function, and $c:\Rbb^n\to\Rbb^m$ is the constraint functions.

One of the most successful and widely adopted CBO methods is the expected constrained improvement (ECI)~\cite{schonlau1998global,gelbart2014bayesian,gardner2014}, where the expected improvement of the objective
is multiplied by the probability of constraints being satisfied. An important reason to its success is that no new statistics routines are needed in its implementation compared to expected improvement, as the probability of feasibility only requires the cumulative 
distribution function (CDF) and probability density function (PDF) of the standard normal distribution. 
Its extension to noisy setting and scalable methods can be found in~\cite{letham2019constrained,eriksson2021scalable}, respectively. 
To account for useful information provided by infeasible samples, integrated conditional expected improvement (IECI)~\cite{bernardo2011}, and expected volume reduction (EVR)~\cite{picheny2014stepwise} are proposed. 
Other than expected improvement, information-based acquisition functions have seen success in the constrained setting with methods such as the predictive entropy search with constraints (PESC)~\cite{hernandez2015predictive}. 
In~\cite{zhang2023constrained}, the authors propose an active learning algorithm and a corresponding acqusition function to search for feasible region. 

Augmented Lagrangian (AL) methods are another important group of CBO methods that draw ideas from the AL methods in nonlinear optimization~\cite{gramacy2016modeling}. These methods construct augmented Lagrangian functions based on the values of both the objective and the constraints. The square of a Gaussian process leads to a new distribution that could require approximation by Monte Carlo integration. In the latest iteration of AL methods, Slack-AL~\cite{picheny2016bayesian} and ADMMBO~\cite{ariafar2019admmbo}, a closed form acqusition function has been achieved, while the algorithms themselves have become more complex.
Recently, constraint violation has been proposed as an improvement on ECI algorithms when no feasible samples are available~\cite{jiao2019complete}.
The authors design a separate branch of the algorithm using the constraint functions alone when no feasible points are available. 
In~\cite{zhao2024bayesian}, the authors opted for a smoothed exact penalty function  and choose scaled expected improvement (SEI) as the acquisition function. However, SEI can encounter zero regions, in which case the algorithm resorts to a mean value reduction scheme.

Despite the recent progress, there remain several limitations for majority of the CBO methods. First, many methods, including ECI, IECI  and EVR, require feasible starting points, which is not available for many engineering problems. 
Second, some of the more recent methods that have superior reported performances such as IECI, EVR, AL Bayesian optimization~\cite{gramacy2016modeling} and PESC  do not have have closed forms, and therefore require computationally expensive estimate routines in their acquisition functions. 
Third, the constraints in methods such as ECI are assumed to be conditionally independent, which can potentially hinder their abilities to find feasible optimal solutions.
Finally, some methods that enjoy closed form acqusition functions and excellent reported performance such as ADMMBO are not readily available in different coding languages as they require additional libraries and more complicated implementation.

In this paper, we propose two novel merit constrained Bayesian optimization (MCBO) algorithms using expected merit improvement (EMI) acqusition functions, both based on the penalty merit function and expected improvement. 
We further propose a unified constrained Bayesian optimization (UCBO) algorithm by combining EMI and ECI for a unified expected constrained improvement (UECI) acquisition function.
The essence of our algorithms is the additive structure from EMI and a relaxation from the double maximum functions on the constraint functions to one. 
Our algorithms can overcome all the limitations mentioned above. That is, it allows infeasible starting point, has closed-form expressions, do not require constraints to be conditionally independent, and finally, can be easily implemented wherever unconstrained Bayesian optimization is. Moreover, the unified algorithm can also been seen as an extension or improvement to ECI.

This paper is organized as follows. In section~\ref{se:bo}, we describe the general Bayesian optimization method with expected improvement acquisition function and its extension to the constrained setting via ECI.
In Section~\ref{se:mcbo}, we motivate and propose our MCBO and UCBO algorithms and their implementation.
Numerical experiments on three synthetic problems and a data-driven ICF design test problem are presented in section~\ref{se:exp} that illustrate and compare the capabilities of the proposed methods.
Conclusions are given in section~\ref{se:conclusion}.

\section{Bayesian optimization}\label{se:bo}
In Bayesian optimization, a subproblem of~\eqref{eqn:opt-prob} is solved at each iteration, where a surrogate function replaces $f$. 
A surrogate function is an approximation of $f$ constructed based on a collection of available data samples. It is formed and updated at each iteration to obtain a more accurate approximation as the optimization progresses.   
A Gaussian process surrogate model assumes that at sample point $x$, a multivariate jointly Gaussian distribution describes the existing samples  and the design objective~$f$ as follows: 
\begin{equation} \label{eqn:GP-1}
 \centering
  \begin{aligned}
  \begin{bmatrix}
   f(x_1) \\
   \vdots\\
           f(x_T)
  \end{bmatrix} \sim \mathcal{N}\left( 
      \begin{bmatrix}    
                 m(x_1)\\
 \vdots\\
 m(x_T)
      \end{bmatrix},
              \begin{bmatrix}
               k(x_1,x_1) \dots k(x_1,x_T)\\
      \vdots\\
               k(x_T,x_1) \dots k(x_T,x_T)
      \end{bmatrix}\ 
      \right). 
\end{aligned}
\end{equation}
Here, $\mathcal{N}$ denotes the normal distribution and $T$ is the number of samples. The input design parameters or samples are denoted as $x_1,\ldots,x_T$. 
The function $m:\Rbb^N\to\Rbb$ is the user-defined  mean function, often chosen to be $0$ constant function.  
The covariance functions are $k:\Rbb^n\times\Rbb^n\to\Rbb$. 
The posterior Gaussian probability distribution of a new sample point~$x$ can be inferred using Bayes' 
rule:
\begin{equation} \label{eqn:GP-post}
 \centering
  \begin{aligned}
  &f(x) | x, x_{1:T},f(x_{1:T}) \sim \mathcal{N} (\mu(x),\sigma^2(x))\\
  &\mu(x)\ =\ k(x,x_{1:T}) k(x_{1:T},x_{1:T})^{-1} \left(f(x_{1:T})-m(x_{1:T}) \right) + m(x)\\
  &\sigma^2(x)\ =\
k(x,x)-k(x,x_{1:T})k(x_{1:T},x_{1:T})^{-1}\sigma(x_{1:T},x)\ ,
\end{aligned}
\end{equation}
where $x_{1:T}=[x_1,\dots,x_T]$ and 
$k(x_{1:T},x_{1:T}) = [k(x_1),\dots,k(x_T); \dots; k(x_T,x_1),\dots,k(x_T,x_T)]$. 
The function $\mu$ and $\sigma^2$ are the posterior mean and variance, respectively.
The covariance functions, or the kernels, of the Gaussian process impact the accuracy of the surrogate model significantly. 
Choices of the functions include the squared exponential covariance function, matern functions, etc.
The hyperparameters of the kernel, often denoted as $\theta$, are obtained through the fitting or training process of the Gaussian process. 

Selection of a suitable acquisition function is critical to an efficient and effective Bayesian optimization algorithm
as it determines the location of the new samples.
One of the most popular acquisition functions is expected improvement (EI)~\cite{brochu2010}, which is computed from the conditional expectation $\Ebb$ of an improvement function defined as
\begin{equation} \label{eqn:improvement}
 \centering
  \begin{aligned}
    I(x) = \max\{ f(x^+) -f(x),0  \}, 
  \end{aligned}
\end{equation}
where 
\begin{equation} \label{eqn:improvement-2}
 \centering
  \begin{aligned}
     x^+= \underset{\substack{x_i\in x_{1:T}}}{\text{argmin}} f(x_i).
  \end{aligned}
\end{equation}
In other words, $x^+$ has the best objective value among existing samples.
The expected improvement is the expectation of~\eqref{eqn:improvement} conditioned on existing samples.
An important advantage of the expected improvement function is that is has a closed form:
\begin{equation} \label{eqn:EI-1}
 \centering
  \begin{aligned}
       EI(x) = \begin{cases}
       (f(x^+)-\mu(x))\Phi(Z)+\sigma(x)\phi(Z), & \sigma(x)>0,\\
       0, & \sigma(x) = 0,    
                 \end{cases}
  \end{aligned}
\end{equation}
where 
\begin{equation} \label{eqn:EI-2}
 \centering
  \begin{aligned}
       Z\ =\ \begin{cases}
       \frac{f(x^+)- \mu(x)}{\sigma(x)}, & \ \sigma(x)>0\\
       0, & \ \sigma(x) = 0\ .    
                 \end{cases}
  \end{aligned}
\end{equation}
A trade-off parameter~$\xi\geq 0$ can be added to~\eqref{eqn:EI-1} for more control over global search and local optimization~\cite{lizotte2008,jones2001taxonomy}. 
The functions 
$\phi$ and~$\Phi$ are the PDF and CDF of the standard normal distribution, respectively.
Therefore, they are readily available in most of the coding languages that one might use to implement Bayesian optimization.  

The next sample point $x_k$ in the Bayesian optimization algorithm is chosen by maximizing the acquisition function over $C$. 
At the $k$th iteration, this means 
\begin{equation} \label{eqn:acquisition-1}
 \centering
  \begin{aligned}
      x_{k} = \underset{\substack{x_i\in C}}{\text{argmin}}  EI_k (x),
  \end{aligned}
\end{equation}
where $EI_k$ is the expected improvement function using samples accumulated up to the $k$th iteration.
In order to solve~\eqref{eqn:acquisition-1}, optimization algorithms including L-BFGS or random 
search can be used. 

Expected constrained improvement (ECI)~\cite{bernardo2011,gelbart2014bayesian,gardner2014} is one of the most widely adopted acqusition functions for constrained optimization problems due to its effectiveness and ease to implement. The constrained improvement function is defined as  
\begin{equation} \label{eqn:c-improvement}
 \centering
  \begin{aligned}
    I_C(x) = \Delta(x)\max\{ f(x^+) -f(x),0  \} = \Delta(x) I(x), 
  \end{aligned}
\end{equation}
where $\Delta(x)\in\{0,1\}$ is a feasibility indicator function that has the value $1$ if $x$ is a feasible point and $0$ otherwise.
The optimal existing sample $x^+$ is now  
\begin{equation} \label{eqn:c-improvement-2}
 \centering
  \begin{aligned}
     x^+=  \underset{\substack{x_i\in x_{1:T}\\c(x)\geq 0}}{\text{argmin}} f(x_i).
  \end{aligned}
\end{equation}
It is clear from~\eqref{eqn:c-improvement-2} that a feasible point in the initial data set is required.
The expected value of the constrained improvement function has a simple closed form
\begin{equation} \label{eqn:cEI-1}
 \centering
  \begin{aligned}
       EI_C(x) = PF(x)EI(x),
  \end{aligned}
\end{equation}
where $PF(\cdot)$ is the probability of a feasible $x$ and it follows a simple univariate Gaussian CDF. 
For multiple inequality constraints, \textit{i.e.}, $m>1$, they are often assumed to be conditionally independent and the probability of $x$ being feasible is computed as 
\begin{equation} \label{eqn:cEI-2}
 \centering
  \begin{aligned}
     p(c(x)\geq 0) = \prod p(c_i(x)\geq 0), 
  \end{aligned}
\end{equation}
where $c_i:\Rbb^n\to\Rbb$ is the $i$th component of the constraint functions.
\section{Merit constrained Bayesian optimization}\label{se:mcbo}
We introduce our proposed acqusition functions and algorithms in this section. 
Our idea is to construct the acquisition function based on merit functions that have been widely adopted in modern nonlinear optimization methods such as sequential quadratic programming (SQP)~\cite{optimization}. While in theory all merit functions can be applied here, we focus on the penalty merit function in this paper. 

\subsection{Expected merit improvement}\label{se:mei}
The penalty merit function can be defined as
\begin{equation} \label{eqn:merit-function}
 \centering
  \begin{aligned}
    \varphi(x) = f(x) + \alpha c^+(x) := f(x) + \alpha \max\{-c(x), 0\}, 
  \end{aligned}
\end{equation}
where $\alpha\geq 0$ is a penalty parameter. It is generally well-understood how $\alpha$ should be updated~\cite{optimization}.
Applying~\eqref{eqn:merit-function} to define the improvement function has been considered previously~\cite{zhao2024bayesian,jiao2019complete}. However, as is evident in its form, this would lead to a maximum of non-Gaussian distribution in the form of 
\begin{equation} \label{eqn:m-improvement-0}
 \centering
  \begin{aligned}
    I_{m,0}(x) = \max\{ f(x^+) + \alpha c^+(x^+) - [f(x) + \alpha c^+(x)], 0 \}, 
  \end{aligned}
\end{equation}
where 
\begin{equation} \label{eqn:m-improvement-xp}
 \centering
  \begin{aligned}
     x^+= \underset{\substack{x_i\in x_{1:T}}}{\text{argmin}} f(x_i)+\alpha c^+(x_i) = \underset{\substack{x_i\in x_{1:T}}}{\text{argmin}} \varphi(x_i).
  \end{aligned}
\end{equation}
The expected value of~\eqref{eqn:m-improvement-0} no longer has a closed form and requires expensive sampling to estimate. 
To resolve this issue, we look back at~\cite{jones1998efficient}, where the authors argued for the expected value of an improvement function, \textit{i.e.}, EI, instead of the expected value of the objective itself for the much superior globalization property of the former. However, for the constraint violation functions $c^+$, a maximum function similar to that in the improvement function is already present. Hence, the global search for $c$ could also exist in the expectation of $c^+$. Therefore, we propose to relax form~\eqref{eqn:m-improvement-0} into two merit improvement functions. The first form is 
\begin{equation} \label{eqn:m-improvement-1}
 \centering
  \begin{aligned}
    I_{m,1}(x) =& \max\{ f(x^+)-f(x),0\} + \alpha c^+(x^+) - \alpha c^+(x)\\ 
               =& \max\{ f(x^+)-f(x),0\} + \alpha \max\{-c(x^+),0\} - \alpha \max\{-c(x), 0 \}. 
  \end{aligned}
\end{equation}
And the second form is 
\begin{equation} \label{eqn:m-improvement-2}
 \centering
  \begin{aligned}
    I_{m,2}(x) =& f(x^+)-f(x) + \alpha c^+(x) - \alpha c^+(x)\\
               =& \varphi(x^+) -f(x) - \alpha \max\{-c(x), 0 \}. 
  \end{aligned}
\end{equation}
The essential idea of~\eqref{eqn:m-improvement-1} and~\eqref{eqn:m-improvement-2} is that since the expected values of the constraint violation functions already use the uncertainty information from $c$ via the maximum functions, we might be able to maintain good globalization property of the acqusition function while removing the additional maximum function in~\eqref{eqn:m-improvement-0}.
We point out that the second form~\eqref{eqn:m-improvement-2} does not use uncertainty information from the objective, thus could be less effective if the objective functions have many local minimum.

We call the expected value of~\eqref{eqn:m-improvement-1} the form 1 expected merit improvement (EMI), denoted as $EI_{M,1}$, which has the closed form
\begin{equation} \label{eqn:m-ei-1}
 \centering
  \begin{aligned}
  EI_{M,1}(x) =& (f(x^+)-\mu_f(x))\Phi(Z_f) +\sigma_f \phi(Z_f) + \alpha c^+(x^+) - \alpha [-\mu_{c}(x)\Phi(Z_c) +\sigma_{c} \phi(Z_c)],\\ 
    =& (f(x^+)-\mu_f(x))\Phi(Z_f) +\sigma_f \phi(Z_f) + \alpha c^+(x^+) + \alpha \mu_{c}(x)\Phi(Z_c) -\alpha\sigma_{c} \phi(Z_c),\\ 
  \end{aligned}
\end{equation}
where 
\begin{equation} \label{eqn:m-ei-2}
 \centering
  \begin{aligned}
       Z_f =  \frac{f(x^+)-\mu_f(x)}{\sigma_f}, \ Z_c = \frac{-\mu_{c}(x)}{\sigma_{c}}. 
  \end{aligned}
\end{equation}
The expression for  form 2 EMI of~\eqref{eqn:m-improvement-2} is 
\begin{equation} \label{eqn:m-ei-3}
 \centering
  \begin{aligned}
  EI_{M,2}(x)  =& \varphi(x^+)-\mu_f(x) - \alpha [-\mu_{c}(x)\Phi(Z_c) +\sigma_{c} \phi(Z_c)]\\ 
                =&\varphi(x^+)-\mu_f(x) + \alpha \mu_{c}(x)\Phi(Z_c) -\alpha\sigma_{c} \phi(Z_c),
  \end{aligned}
\end{equation}
where $Z_c$ is the same as in~\eqref{eqn:m-ei-2}.
Both acqusition functions can be easily extended to multiple constraints. The form 1 EMI becomes 
\begin{equation} \label{eqn:m-ei-4}
 \centering
  \begin{aligned}
   EI_{M,1}(x) &= (f(x^+)-\mu_f(x))\Phi(Z_f) +\sigma_f \phi(Z_f) + \alpha \sum_{j=1}^m c^+_j(x^+) \\
                      &+ \alpha \sum_{j=1}^m [\mu_{c_j}(x)\Phi(Z_{c_j}) -\sigma_{c_j} \phi(Z_{c_j})],\\ 
  \end{aligned}
\end{equation}
where for $j=1,\dots,m$,
\begin{equation} \label{eqn:m-ei-5}
 \centering
  \begin{aligned}
        Z_{c_j} = \frac{-\mu_{c_j}(x)}{\sigma_{c_j}}. 
  \end{aligned}
\end{equation}
Similarly, the second form is 
\begin{equation} \label{eqn:m-ei-6}
 \centering
  \begin{aligned}
   EI_{M,2}(x) =& \varphi(x^+)-\mu_f(x)+  \alpha \sum_{j=1}^m [\mu_{c_j}(x)\Phi(Z_{c_j}) -\sigma_{c_j} \phi(Z_{c_j})].\\ 
  \end{aligned}
\end{equation}
The penalty merit function~\eqref{eqn:merit-function} with multiple constraints can be defined as
\begin{equation} \label{eqn:merit-function-2}
 \centering
  \begin{aligned}
    \varphi(x) = f(x) + \alpha \sum_{j=1}^m \max\{-c_j(x), 0\}. 
  \end{aligned}
\end{equation}
The optimal existing sample is computed using the same formula~\eqref{eqn:m-improvement-xp} with the updated merit function~\eqref{eqn:merit-function-2}.
The MCBO algorithms are given in Algorithm~\ref{alg:mcbo}.
\begin{algorithm}[H]
 \caption{Merit constrained Bayesian optimization}\label{alg:mcbo}
  \begin{algorithmic}[1]
	  \State{Choose parameter $\alpha\geq 0$ and its update rule. Choose either~\eqref{eqn:m-ei-4} or~\eqref{eqn:m-ei-6} as the acquisition function. Choose $N_0$ initial samples $x_i,~i\in[0,N_0]$. Compute $f(x_i)$ and $c(x_i)$.} 
	  \State{Train Gaussian Process surrogate models for the objective and constraint functions on the initial samples.}
  \For{$k=0,1,2,\dots$}
	  \State{Evaluate the acquisition function and find the new sample point based on~\eqref{eqn:acquisition-1}.}
	  \State{Run (numerical) experiments with the new sample $x_k$ if necessary.}
	  \State{Evaluate the objective $f(x_{k})$ and constraint functions $c(x_k)$. Update $\alpha$ if needed. \;}
	  \State{Retrain the surrogate models with the addition of the new sample $x_k$, $f(x_k)$ and $c(x_k)$.\;}
	  \State{Solve the optimization problem~\eqref{eqn:opt-prob} with the updated surrogate model.
	  \If {Stopping criteria satisfied} Exit
	  } 
  \EndFor
  \end{algorithmic}
\end{algorithm}
We include line $5$ in the algorithm because new numerical simulations might need to be run with $x_k$, as in many engineering applications.
Line $8$ of Algorithm~\ref{alg:mcbo} can be solved through exhaustive search on a random, sufficiently large set of sample points in $C$ or gradient-based methods such as L-BFGS. 
The specific stopping criteria and the update rule for $\alpha$ are left to be chosen by users to allow maximum flexibility of the algorithms. 
We remark that the algorithmic parameter $\alpha$ can be changed, \textit{e.g.}, based on whether the constraints are satisfied, throughout the optimization process. Additionally, each constraint can have a different $\alpha$ for better performance.

\subsection{A unified algorithm with with ECI}
For the penalty merit function to produce equivalent solution to the original problem~\eqref{eqn:opt-prob}, the parameter $\alpha$ needs to be large enough~\cite[Theorem 17.3]{optimization}. 
As a result, the additive structure of the merit function in~\eqref{eqn:m-ei-4} and~\eqref{eqn:m-ei-6} might lead to more exploration within the region that current constraint surrogate models expect to be feasible.
If the boundary of the feasible regions are complex, the algorithms could become less efficient.
Meanwhile, the multiplication structure in the ECI method, without the penalty $\alpha$, can lead to more exploration where the current constraint surrogate models deem close to the boundary but infeasible. 
Therefore, we propose a unified acqusition function that combines ECI with EMI through algorithmic parameters, referred to as the unified expected constrained improvement (UECI). Its improvement function is given by  
\begin{equation} \label{eqn:u-improvement-1}
 \centering
  \begin{aligned}
    I_{u}(x)  =& (1-\beta) \Delta(x)\max\{ f(x^+)-f(x),0\}  \\
               &+ \beta \left[ \max\{ f(x^+)-f(x),0\} + \alpha\sum_{j=1}^m c_j^+(x^+) - \alpha \sum_{j=1}^m c^+(x)\right],\\ 
  \end{aligned}
\end{equation}
where $\beta\in[0,1]$ is a parameter. When $\beta=0$, we recover ECI~\eqref{eqn:c-improvement}, and when $\beta=0$,~\eqref{eqn:u-improvement-1} reduces to~\eqref{eqn:m-improvement-1}. 
We note that many other combinations are possible using, \textit{e.g.},~\eqref{eqn:m-ei-5} and other algorithmic parameter design.
The UECI expression is  
\begin{equation} \label{eqn:ueci-1}
 \centering
  \begin{aligned}
    EI_{UC}  =& (1-\beta) PF(x)EI(x) + \beta EI_{m,1}(x),\\ 
  \end{aligned}
\end{equation}
based on~\eqref{eqn:cEI-1},~\eqref{eqn:cEI-2}, and~\eqref{eqn:m-ei-4}.
The UECI function continues to be closed-formed, making its implementation straightforward.
In this paper, we devise a simple update rule for $\beta$ that relies on the number of feasible points.
Specifically, let $n_f$ be the number of existing feasible samples and $N_f$ the threshold.
A unified constrained Bayesian optimization (UCBO) algorithm based on UECI~\eqref{eqn:ueci-1} is given in Algorithm~\ref{alg:uecibo} 
\begin{algorithm}[H]
 \caption{Unified constrained Bayesian optimization}\label{alg:uecibo}
  \begin{algorithmic}[1]
	  \State{Choose parameter $\alpha\geq 0$ and $N_f\geq 1$. Set $\beta=1$. Choose $N_0$ initial samples $x_i,~i\in[0,N_0]$. Compute $f(x_i)$ and $c(x_i)$.} 
	  \State{Train Gaussian Process surrogate models for the objective and constraint functions on the initial samples.}
  \For{$k=0,1,2,\dots$}
          \If {$n_f \geq N_f$} 
            \State{Set $\beta=0$.}
          \EndIf
	  \State{Evaluate the acquisition function~\eqref{eqn:ueci-1} and find the new sample point with~\eqref{eqn:acquisition-1}.}
	  \State{Run (numerical) experiments with the new sample $x_k$ if necessary.}
	  \State{Evaluate the objective $f(x_{k})$ and constraint functions $c(x_k)$. Update $\alpha$ if needed. \;}
	  \State{Retrain the surrogate models with the addition of the new sample $x_k$, $f(x_k)$ and $c(x_k)$.\;}
	  \State{Solve the optimization problem~\eqref{eqn:opt-prob} with the updated surrogate model.}
	  \If {Stopping criteria satisfied} 
              \State{Exit}
	  \EndIf
  \EndFor
  \end{algorithmic}
\end{algorithm}
Algorithm~\ref{alg:uecibo} offers a simplistic switch between ECI and MCBO methods. The choices for $\beta$ do not need to be binary.
We leave more complex designs for update rules of both $\alpha$ and $\beta$, as well as algorithms, for future work.
Algorithm~\ref{alg:mcbo} and~\ref{alg:uecibo} have been implemented in Python using the \texttt{scikit-learn} machine learning library~\cite{scikit-learn,sklearn_api} and MATLAB.  

\section{Experiments}\label{se:exp}
In this section, we demonstrate and evaluate the performance of MCBO and UCBO algorithms on three synthetic problems that have been studied and compared before~\cite{gardner2014,gramacy2016modeling,picheny2016bayesian,ariafar2019admmbo}. Additionally, we apply MCBO and UCBO algorithms to a data-driven ICF design problem based on complex finite element simulation data and practical design objective and constraints. 

For the synthetic problems, we compare our proposed algorithms with the state-of-the-art constrained Bayesian optimization methods ADMMBO method and the ECI~\cite{gardner2014} method with additional random sampling. As reported in~\cite{ariafar2019admmbo}, ADMMBO is considered to be one of the best performing CBO algorithms, as it is built upon previous AL methods.    
Meanwhile, ECI is one of the most widely adopted CBO methods thanks to its simplicity and  effectiveness. 
For a fair comparison, the ADMMBO results in this paper is obtained using the open-source code the authors of~\cite{ariafar2019admmbo} published, where the ADMMBO parameters and initial data sets are unchanged to maximize its performance. 
Therefore, readers can find the detailed parameter choices and their discussion for ADMMBO in~\cite{ariafar2019admmbo}. 
If no initial samples are feasible, we allow the ECI implementation to call Latin hypercube sampling to generate new samples until a feasible one is found. 
The MCBO, UCBO, and ECI algorithms are implemented in both Python and MATLAB, while ADMMBO is currently available in MATLAB.

Each method is run $100$ times to obtain a more comprehensive comparison of performances.
ADMMBO has two initial samples for the objective and constraints, respectively, while four initial Latin hypercube samples are used for other algorithms. 
The median of the best feasible objectives among the $100$ runs is reported against the number of iterations.
 Additionally, we report the $25$ and $75$ percentiles of the best feasible objectives from our proposed algorithms to illustrate their robustness. 
For the ADMMBO method, the iteration number refers to its inner iteration, which is given a budget of $100$ by the authors of~\cite{ariafar2019admmbo}.
Thanks to a stopping criterion, ADMMBO sometimes converges before the budget is used up.
The budgets of other algorithms are chosen based on when the best feasible objectives stabilize, as more iterations do not change their performances.

We choose to compare the number of iterations since the evaluation of objective and constraint functions can be done independently. More importantly, in many engineering applications, \textit{e.g.}, ICF design optimization, a new sample in the design space requires a new complex and computationally expensive numerical simulation, while the design objective and constraints evaluations are analytical and thus much less costly. We do point out that ADMMBO has the advantage that only one function, either the objective or the constraint, is evaluated at one inner iteration.    
\subsection{Example with small feasible region}
\begin{figure}
  \centering
   \includegraphics[width=\textwidth]{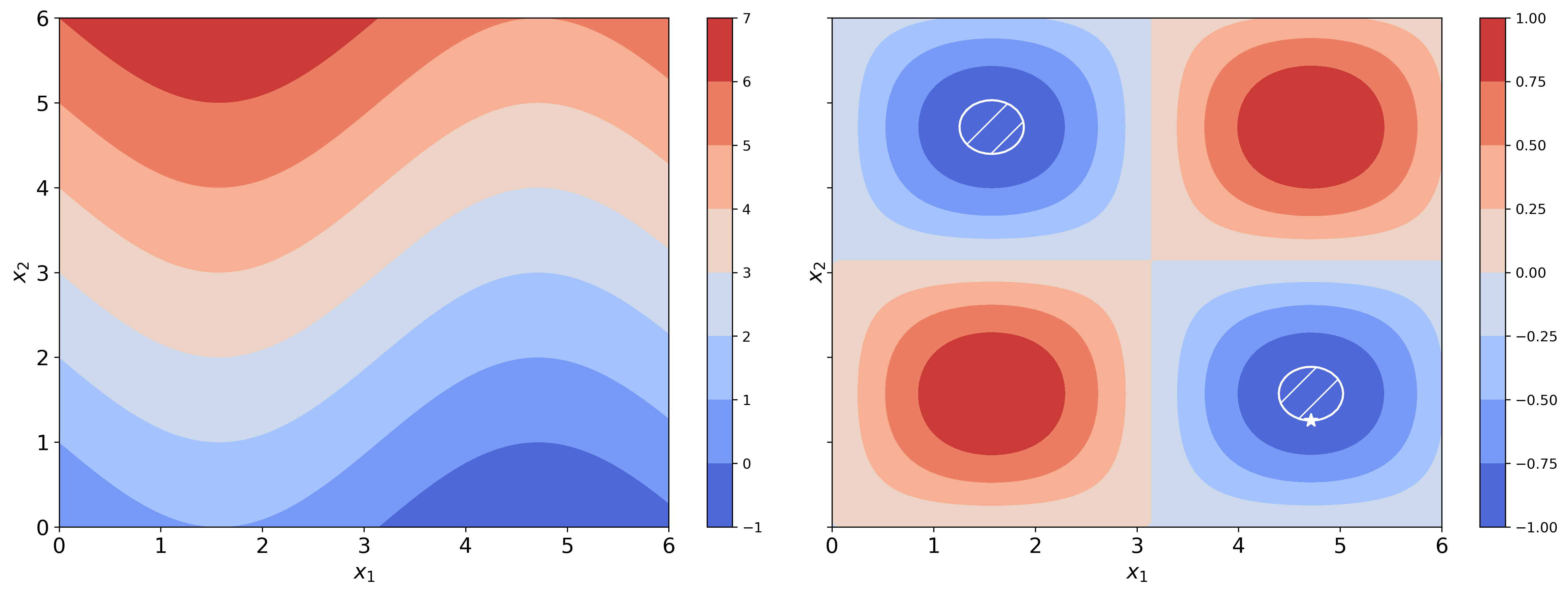}
  \caption{Contour plots for the objective function (left) and constraint function (right) for example 1. 
            The feasible region is marked on the constraint plot, along with the global optimum ($*$ sign).}
\label{fig:ex1contour}
\end{figure}

The mathematical expression for the first example is 
\begin{equation} \label{eqn:ex1}
 \centering
  \begin{aligned}
	  &\underset{\substack{x}\in C}{\text{minimize}} 
	  & & \sin(x_1)+ x_2, \\
          &\text{subject to}
	  & & -\sin(x_1)\sin(x_2) -0.95 \geq 0,
  \end{aligned}
\end{equation}
where $C=[0,6]\times[0,6]$. This is a classic problem that has been used in literature due to its small feasible region and the existence of distinct local and global minimum. To demonstrate this point, the feasible region and the objective contour are plotted in figure~\ref{fig:ex1contour}. The local minimum has a value of $5.4$ and the global minimum has the value of $0.25$, obtained through exhaustive random sampling of $10000$ points.

Given the small feasible region, we choose $\alpha=20$ for form 1~\eqref{eqn:m-ei-4} of MCBO, $\alpha=5$ for form 2~\eqref{eqn:m-ei-6}, and $\alpha=20, N_f=2$ for UCBO. MCBO and UCBO algorithms are given $60$ iterations, while ECI and ADMMBO are given $100$. 

\begin{figure}
  \centering
   \includegraphics[width=0.9\textwidth]{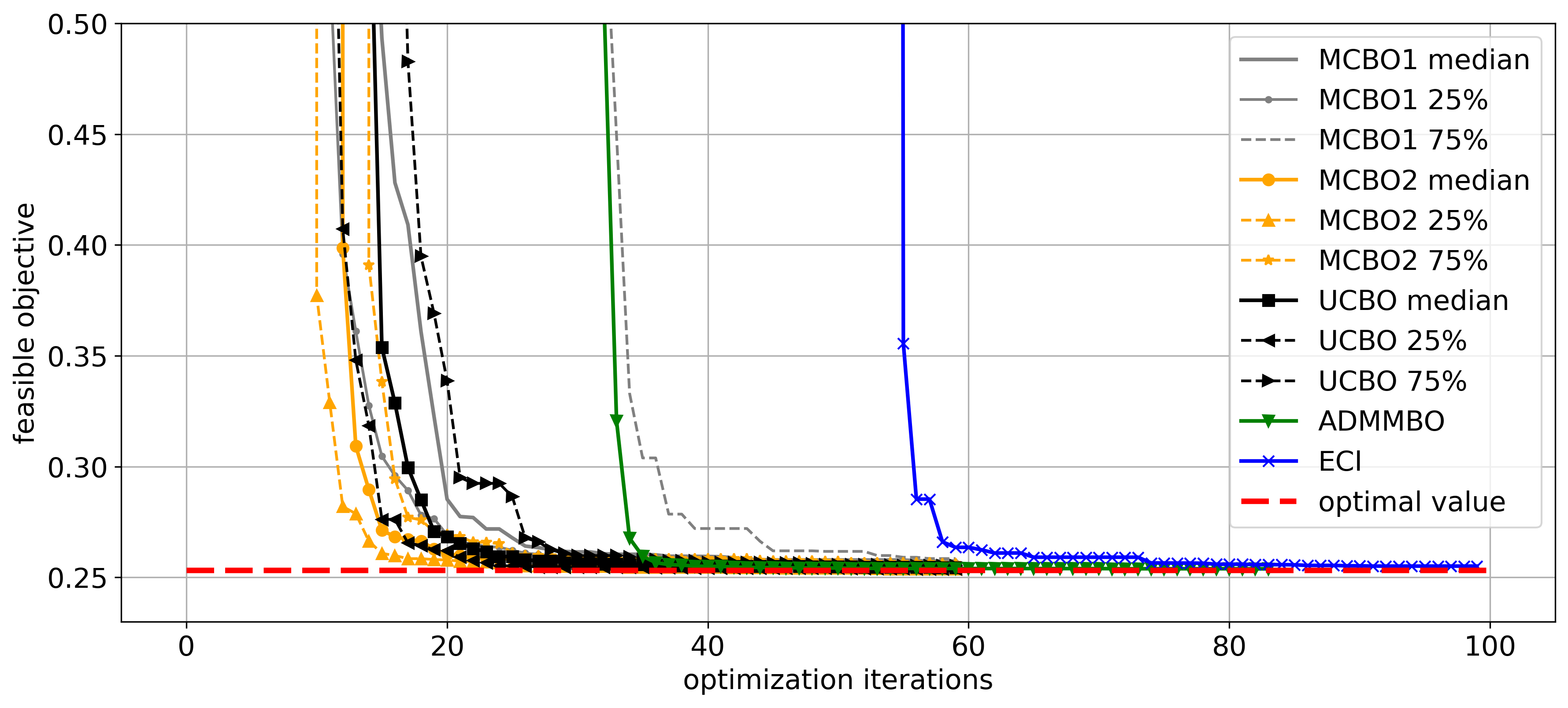}
  \caption{Optimization result of MCBO, UCBO, ECI (median) and ADMMBO (median) algorithms for example 1, shown in figure~\ref{fig:ex1contour}. 
           The best feasible objective is plotted against the number of iterations, or equivalently, number of samples.}
  \label{fig:ex1result}
\end{figure}

Figure~\ref{fig:ex1result} shows the median of best feasible objective against number of iterations among $100$ runs. Form 2 of MCBO and UCBO converge the fastest to the optimal value, followed closely by form 1 of MCBO,  while ECI is the slowest. This is likely due to the small feasible region that ECI is incapable of finding effectively. All algorithms succeed in escaping the local minimum and searching in the neighborhood of the global optimum. 
The $25$ and $75$ percentiles of the best feasible objectives from MCBO and UCBO algorithms over the $100$ runs illustrate the robustness and reliability of the proposed algorithms.
It is worth pointing out that a more sophisticated update rule for $\alpha$ and $\beta$ is likely to improve performance further for the proposed algorithms, a topic for future work. 
\subsection{Example with multiple constraints}
\begin{figure}
  \centering
   \includegraphics[width=\textwidth]{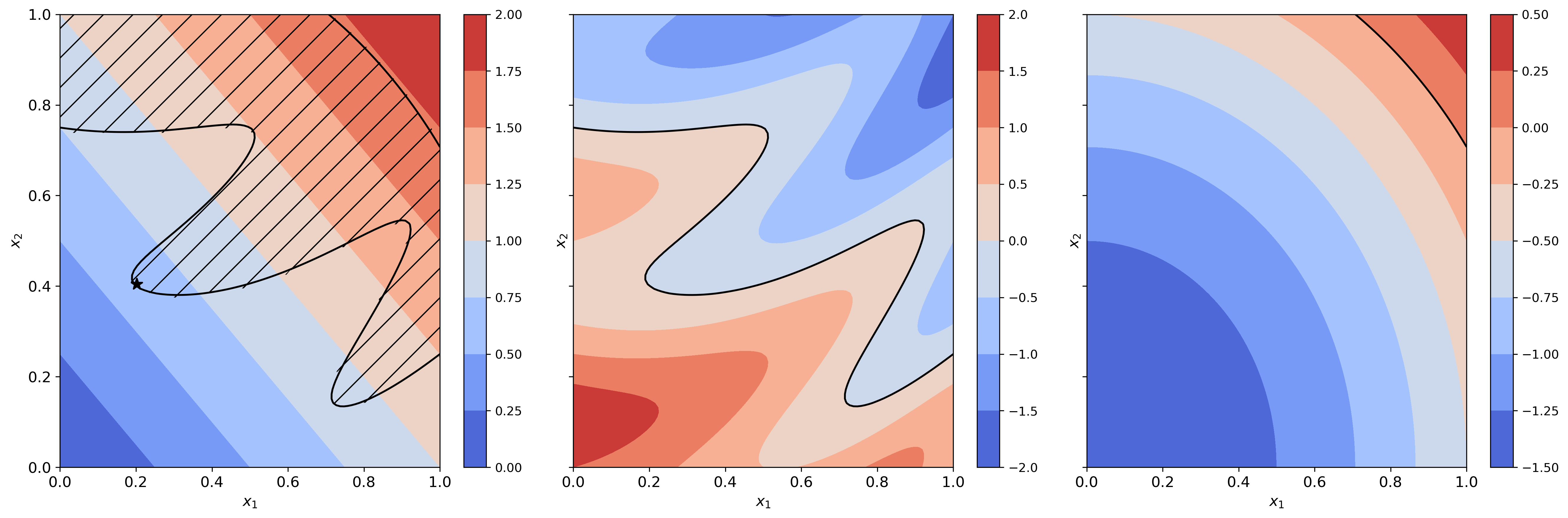}
  \caption{Contour plots for the objective function (left) and the two constraint functions (middle and right) for example 2. 
            The feasible region is marked with black line on the objective contour, along with the global optimum ($*$ sign).
            The contour curve $c(x)=0$ is in black on each constraint contour plot as well.}
\label{fig:ex2contour}
\end{figure}

Next, we consider the test problem 
\begin{equation} \label{eqn:ex2}
 \centering
  \begin{aligned}
	  &\underset{\substack{x}\in C}{\text{minimize}} 
	  & & x_1+ x_2, \\
          &\text{subject to}
	  & & 0.5 \sin( 2\pi (x_1^2-2 x_2))+x_1+2x_2 -1.5 \geq 0,\\
          &&& -x_1^2-x^2_2 +1.5\geq 0,
  \end{aligned}
\end{equation}
where $C=[0,1]\times[0,1]$. The feasible region and the objective contour are plotted in figure~\ref{fig:ex2contour}. 
It is evident from the plot that the direction of a smaller objective is in agreement with the feasible region of the second constraint. Additionally, the feasible region is not as restrictive as the one in example 1.

We choose $\alpha=[2,0.02]$ for form 1 of MCBO, $\alpha=[25,25]$ for form 2 of MCBO, and $\alpha=[100,0.1], N_f=1$ for UCBO.
ADMMBO is allowed the same $100$ inner iterations without modifying the original code to ensure fair comparison.
Based on their performances, form 1 of MCBO, UCBO, and ECI are run for $60$ iterations while form 2 has a budget of $80$ iterations that matches the iterations reported by ADMMBO upon exit. 

\begin{figure}
  \centering
   \includegraphics[width=0.9\textwidth]{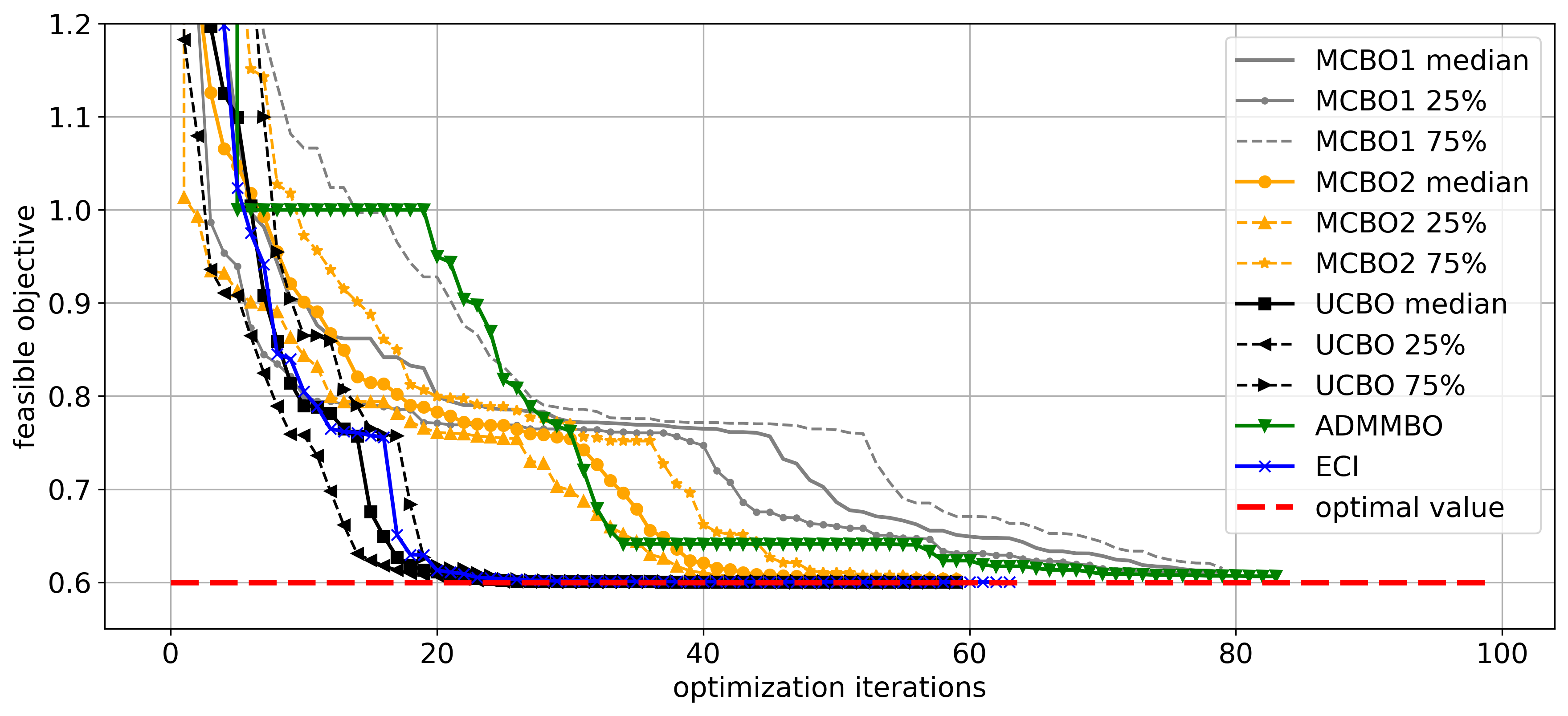}
  \caption{Optimization result of MCBO, UCBO, ECI (median), and ADMMBO (median) algorithms for example 2, shown in figure~\ref{fig:ex2contour}. 
           The best feasible objective is plotted against the number of iterations.}
  \label{fig:ex2result}
\end{figure}

The median of the best feasible objectives are plotted in figure~\ref{fig:ex2result}, together with the best feasible objectives at $25$ and $75$ percentiles for MCBO and UCBO algorithms. Form 2 MCBO finds feasible samples and converges to optimal solution faster than ADMMBO, while form 2 initial outperforms ADMMBO but eventually find the optimum in similar number of iterations. The $25\%$ and $75\%$ curves demonstrate the robustness and stability of the two proposed algorithms.
UCBO is the best performing algorithm among all those tested. 

Perhaps the surprising result is the excellent performance of ECI for this problem in our experiments, contrary to findings in previous papers~\cite{ariafar2019admmbo,gramacy2016modeling,hernandez2015predictive}. For example, in~\cite{ariafar2019admmbo}, ECI does not have feasible samples for the $100$ runs until $40$ iterations, a significantly worse performance than ours.  
In our experiments, given that the feasible region is relative large in this low-dimensional problem, ECI, with the help of continued Latin hypercube sampling, is able to have a feasible sample within the first eight samples. 
It is then capable of finding global minimum very efficiently. 
Since it is unclear precisely how other authors overcome the lack of initial feasible samples for ECI in their experiments, we contemplate that our approach of sampling when there is no feasible sample could be the reason for the improved performance of ECI. 
The UCBO method can fully take advantage the capability of ECI to explore outside of feasible region while working with infeasible initial samples and locating feasible ones quickly. 

\subsection{Four-dimensional problem}
\begin{figure}
  \centering
   \includegraphics[width=0.9\textwidth]{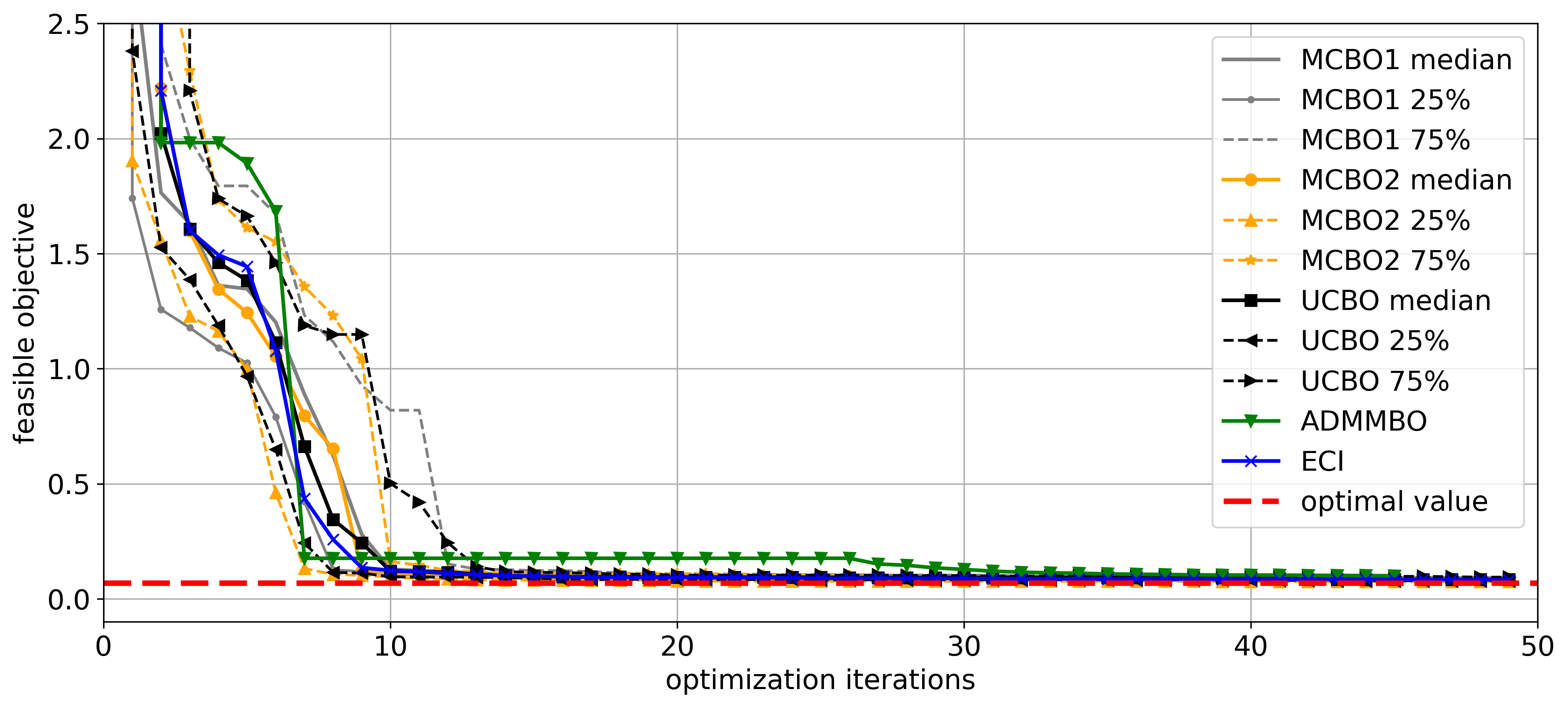}
  \caption{Optimization result of MCBO, UCBO, ECI (median) and ADMMBO (median) algorithms for example 3. 
           The best feasible objective is plotted against the number of iterations.}
  \label{fig:ex3result}
\end{figure}

To examine the performance of our proposed algorithms for problems in higher dimensions, we test a four-dimensional problem that previously appeared in~\cite{picheny2016bayesian,ariafar2019admmbo} 
\begin{equation} \label{eqn:ex4d}
 \centering
  \begin{aligned}
	  &\underset{\substack{x}\in C}{\text{minimize}} 
	  & & \sum_{i=1}^4 x_i, \\
          &\text{subject to}
	  & & \frac{1}{0.8387}\left[  \sum_{i=1}^4 E_i \exp( -\sum_{j=1}^4 A_{ji}(x_j-P_{ji})^2 ) -1.1\right]\geq 0,\\
  \end{aligned}
\end{equation}
where $C=[0,1]^4$ and
\begin{equation} \label{eqn:ex4d-matrix}
 \centering
  \begin{aligned}
        E = \begin{bmatrix*}
             &1\\&1.2\\&3\\&3.2
            \end{bmatrix*},\
        A = \begin{bmatrix*}
             &10 \ &0.05 \ &3 \ &17\\
             & 3 \ &10 \ &3.5 \ &8\\
             &17 \ &17 \ &1.7\ &0.05\\
             &3.5 \ &0.1 \ &10 \ &10
            \end{bmatrix*},\
        P =\begin{bmatrix*}
            &0.131 \ &0.232 \ &0.234 \ &0.404\\
            &0.169 \ &0.413 \ &0.145 \ &0.882\\
            &0.556 \ &0.83 \ &0.352 \ &0.873\\
            &0.012 \ &0.373\ &0.288 \ &0.574 \\
            \end{bmatrix*}.
  \end{aligned}
\end{equation}
We choose an iteration budget of $50$  for all five algorithms.
We let $\alpha=0$ for the first $10$ iterations of form 1 MCBO and $\alpha=0.01$  for the remaining 40 iterations, while we set $\alpha=0.01$ for form 2 MCBO.
Similarly, we let $\alpha=0.01$ and $N_f=2$ for UCBO.
The result of the optimization is shown in figure~\ref{fig:ex3result}, where the same curves as in example 1 and 2 are plotted. 
Among the median results, form 1 MCBO algorithm is the first to find feasible samples. 
Moreover, all algorithms quickly find values close to the global minimum, except for ADMMBO, which might be caused by the constraint tolerance $0.01$ set by the ADMMBO run.
Therefore, the proposed algorithms are capable of top-tier performance in higher dimensions as well. 

It is worth noting that the ADMMBO results for all three synthetic problems here is very similar to those presented in~\cite{ariafar2019admmbo}, where it is the best performing algorithm among all those tested. 
\subsection{ICF design problem}
Finally, we present a data-driven ICF design optimization problem, where the objective and constraint functions are results from complex simulation of multiphysics problems with no analytical expressions. 
The objective and constraint functions are based on simulation data collected using HYDRA, a multi-physics simulation code developed at LLNL~\cite{marinak2001three}.
To perform an efficient study, HYDRA simulation is not called directly by our algorithms.
Rather, a pre-trained neural networks fitted to a HYDRA simulation database of one-dimensional capsule simulations of design variations around NIF shot N210808 is used. It is the first ICF experiment to exceed the Lawson fusion condition and generate more than 1 MJ of energy~\cite{N210808-lawson}. The optimization process would be prohibitively expensive if HYDRA is run natively.

We select design variables `\texttt{t\_2nd}' and `\texttt{sc\_peak}', which are the timing of the second shock and the strength of the peak of the radiation drive. For a figure that illustrates the relationship between the two variables and the radiation temperature, see~\cite{wang2023multifidelity}. The design objective is to maximize the nuclear yield, or equivalently minimizing the negative yield.  
The two constraint functions are chosen to be $c_1$ = `\texttt{adiabat}' and $c_2$= `\texttt{vImp}'.
The implosion velocity vImp (km/sec) is the maximum inward velocity of the fusion fuel achieved during the implosion, and the adiabat is a unit-less measure $>1$ of the entropy of the fuel at the time of peak implosion velocity. In one-dimensional models (like simulated here), implosions with lower values of adiabat and higher values of velocity tend to achieve higher yield; however, more complex models show that additional performance-limiting physics can turn on if either the adiabat becomes too low or the velocity becomes too high. As such, it is often desirable to search for designs that have constraints on these functions.
Therefore, we impose the constraints $c_1(x)\geq 4.25$ and $c_2(x)\leq 350$.

\begin{figure}
  \centering
   \includegraphics[width=\textwidth]{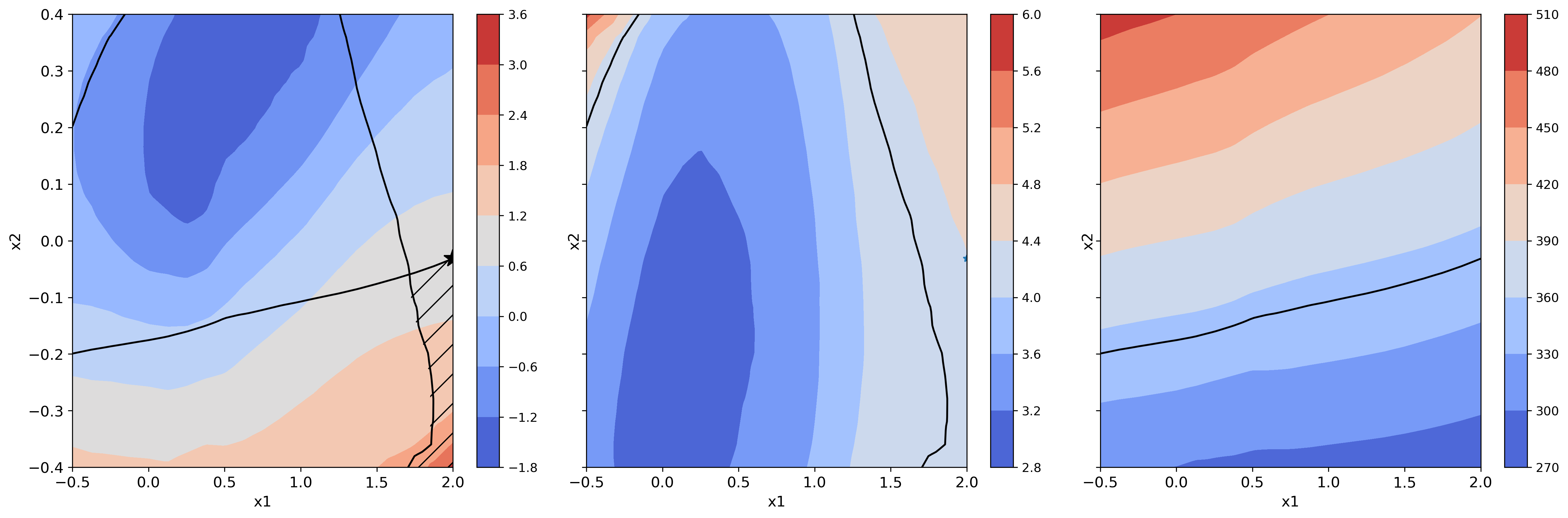}
  \caption{Contour plots for the objective function (left) and the two constraint functions (middle and right) from HYDRA data set. 
            The feasible region is marked with black dashed lines on the objective contour, along with the global optimum ($*$ sign).
            The contour curve $c(x)=0$ is in black curve on each constraint contour plot as well.}
\label{fig:hydracontour}
\end{figure}

The contour plots for the objective (negative yield) and constraints are given in figure~\ref{fig:hydracontour}, along with the feasible region. The optimal feasible objective value is obtained through exhaustive sampling of the two-dimensional design space. 
As evident in figure~\ref{fig:hydracontour}, the feasible region is relatively small for this problem.

\begin{figure}
  \centering
   \includegraphics[width=0.9\textwidth]{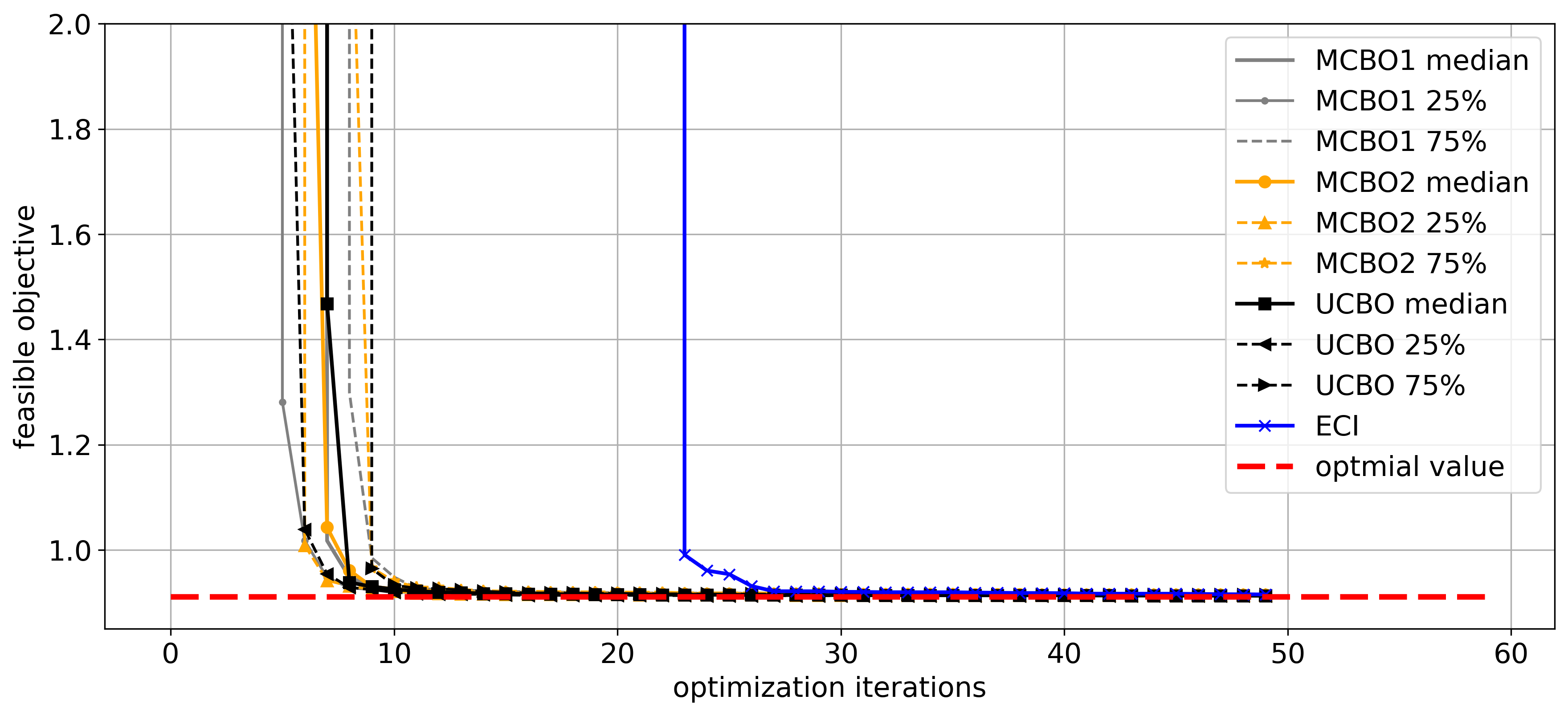}
  \caption{Optimization result of two MCBO, UCBO and ECI (median) algorithms for the ICF design problem in figure~\ref{fig:hydracontour}. 
           The mean/25\%/75\% best feasible objectives of MCBO and UCBO methods are plotted against the number of iterations.}
  \label{fig:hydraresult}
\end{figure}

We adopt to only apply MCBO, UCBO, and ECI algorithms for this problem since they are the most straightforward to implement in an existing and complex code base. For all three proposed algorithms, we choose $\alpha=[1,1]$ and set $N_f=3$ for UCBO due to the small feasible region.
We allow ECI to randomly sample until a feasible point is found. A total budget of $50$ iterations is assigned to all algorithms tested. The result is given in figure~\ref{fig:hydraresult}, where the MCBO and UCBO algorithms behave similarly but clearly outperform ECI thanks to its ability to make progress on constraint violation from the beginning. 
The model ICF design problem presented here can be extended to a full problem with eight or as high as thirty dimensions. 
An eight-dimensional constrained ICF design problems is currently being developed. 

\section{Conclusion}\label{se:conclusion}
In this paper, we propose constrained Bayesian optimization methods with merit expected improvement acquisition functions that can resolve a large number of issues that constrained expected improvement encounter.
Through relaxation on the improvement function of the constraint violation functions, we can maintain closed-form expressions for the new acquisition functions, while making progress at infeasible samples.
The additive structure between the objective and constraint violation functions allows expectation to be taken efficiently on conditionally dependent constraint functions.
Moreover, the proposed algorithms do not require any new routines to implement compared to unconstrained Bayesian optimization.

We further propose a unified expected constrained improvement function and design a constrained Bayesian optimization algorithm that can take advantage of the benefits of expected constraint improvement function if needed. The unified  approach can be seen as a framework where other methods can be integrated. 

Some questions that have not been answered in this paper include the algorithmic parameter update rule for $\alpha$ and $\beta$. 
More numerical experiments in higher dimension would further validate the performance of the proposed algorithms. The authors will work on these topics in the future.  
\begin{ack}
This work was performed under the auspices of the U.S. Department of Energy by Lawrence Livermore National Laboratory under contract DE--AC52--07NA27344 and the LLNL-LDRD Program under Project tracking No. 23-ERD-017.  Release number LLNL-CONF-861915-DRAFT. The data and codes that support the findings of this study are available from the corresponding author upon request.
\end{ack}

\bibliographystyle{unsrtnat}
\bibliography{bibliography}

\end{document}